\documentclass[12pt]{article}       
\usepackage{amsfonts,amsmath,amssymb}
\usepackage[usenames,dvipsnames]{pstricks}
\usepackage{epsfig}
\usepackage{pst-grad}												 
\usepackage{pst-plot} 												

\usepackage{graphicx}

\newtheorem{Def}{Definition}
\newtheorem{Th}{Theorem}
\newtheorem{cor}[Th]{Corollary}

\newcommand{\cn}{{\rm col_g}}

\newenvironment{Pf}
{\noindent\textbf{Proof.}\ \ }{\hfill $\Box$ \medskip}

\title{Bounds for the game coloring number of planar graphs with a specific girth}

\author {Keaitsuda Maneeruk Nakprasit \\ {\small\em Department of Mathematics, Faculty of Science, Khon Kaen University, 40002, Thailand }\\  
{\small\em E-mail address: kmaneeruk@hotmail.com} 
\and Kittikorn Nakprasit \footnote{Corresponding Author}\\ 
{\small\em Department of Mathematics, Faculty of Science, Khon Kaen University, 40002, Thailand }\\
{\small\em E-mail address: kitnak@hotmail.com}}
\date{}

\begin{document}

\maketitle

\begin{center}{\bf Abstract}\end{center} 
Let ${\rm col_g}(G)$ be the game coloring number of a given  graph $G.$ 
Define the game coloring number of  a family of graphs $\mathcal{H}$ 
as ${\rm col_g}(\mathcal{H}) := \max\{{\rm col_g}(G):G \in \mathcal{H}\}.$  
Let $\mathcal{P}_k$ be the family of planar graphs of girth at least $k.$ 
We show that ${\rm col_g}(\mathcal{P}_7) \leq 5.$ 
This  result extends a result about the coloring number  
by Wang and Zhang \cite{WZ11} (${\rm col_g}(\mathcal{P}_8) \leq 5).$  
We also show that these bounds are sharp by constructing 
a graph $G$ where $G \in {\rm col_g}(\mathcal{P}_k) \geq 5$ for each $k \leq 8$ 
such that ${\rm col_g}(G)=5.$   
As a consequence, ${\rm col_g}(\mathcal{P}_k) = 5$ for $k =7,8.$ 
\indent\indent


\section{Introduction} 

Let $G$ be a simple graph with a vertex set $V(G)$ and an edge set $E(G).$ 
The \textit{coloring game} is  a two-person game described as follows. 
Two players, say Alice and Bob, with Alice playing first alternatively 
colors an uncolored vertex in $G$ with the color from the color set $C$ 
so that any two adjacent vertices have distinct colors. 
Alice wins if all vertices are colored. 
The \textit{game chromatic number} of $G,$ denoted by $\chi_g(G),$ is the least 
cardinality of $C$ in which Alice has a strategy to win the game. 
The game chromatic number was formally introduced by Bodlaender \cite{B91}. 

The \textit{marking game} is also a two-person game.  
Two players, say Alice and Bob, with Alice playing first alternatively 
marks an unmarked vertex of $G$ until all vertices are marked. 
Let $b(v)$ be the number of neighbors of $v$ that are marked before $v$ is marked. 
The \textit{game coloring number} of $G,$ denoted by $\cn(G),$ is the least 
$s$ in which Alice has a strategy to obtain $b(v)+1$ to be at most $s$ for each vertex $v.$ 

The game coloring number was formally introduced by Zhu \cite{Z99} as a tool to 
study the game chromatic number. It is easy to see that $\chi_g(G) \leq \cn(G).$ 
The best known upper bounds for game chromatic numbers of graphs in various families   
are obtained from the upper bounds of game coloring numbers. 

Let $\mathcal{H}$ be a family of graphs. 
The game chromatic number and the game coloring number of $\mathcal{H}$ 
are defined as $\chi_G(\mathcal{H}) := \max\{\chi_g(G):G \in \mathcal{H}\}$ 
and $\cn(\mathcal{H}) := \max\{\cn(G):G \in \mathcal{H}\}.$

The game coloring numbers of various families of graphs, especially planar graphs, 
are widely studied.  
Let $\mathcal{F}$ denote the family of forests, 
$\mathcal{I}_k$ denote the family of interval graphs with clique number $k,$ 
$\mathcal{Q}$ denote the family of outerplanar graphs, 
$\mathcal{PT}_k$ denote the family of partial $k$-trees, and 
$\mathcal{P}$ denote the family of planar graphs. 
It is proved by Faigle et al. \cite{FKKT93}  that $\chi_g(\mathcal{F})=\cn(\mathcal{F})=4,$ 
by Faigle et al. \cite{FKKT93} and Kierstead and Yang \cite{KY05} that $\cn(\mathcal{I}_k)=3k-2,$ 
by Guan and Zhu \cite{GZ99} and  Kierstead and Yang \cite{KY05} that $\cn(\mathcal{Q})=7,$ 
and by  Zhu \cite{Z00}  and Wu and Zhu \cite{WZ08} that  $\cn(\mathcal{PT}_k)=3k+2$ for $k\geq 2.$  

Combining a lower bound from \cite{WZ08} and an upper bound from \cite{Z08} gives 
$11 \leq \cn(\mathcal{P}) \leq 17.$ 
Let $\mathcal{P}_k$ be the family of planar graphs of girth at least $k.$ 
It is proved by Sekiguchi \cite{S14} that $\cn(\mathcal{P}_4) \leq 13,$ 
by He et al. \cite{HEtal02} that $\cn(\mathcal{P}_5) \leq 8,$  
by Kleitman \cite{K06} that $\cn(\mathcal{P}_6) \leq 6,$ 
by Wang and Zhang \cite{WZ11} that $\cn(\mathcal{P}_8) \leq 5,$  
and by Borodin et al. \cite{BIKS09} that $\cn(G) \leq 9$ 
if $G$ is a quadrangle-free planar graph.  

In this paper, we show that ${\rm col_g}(\mathcal{P}_7) \leq 5.$ 
This  result extends a result about the coloring number  
by Wang and Zhang \cite{WZ11}.  
We also show that these bounds are sharp by constructing 
a graph $G$ where $G \in {\rm col_g}(\mathcal{P}_k) \geq 5$ for each $k \leq 8$ 
such that ${\rm col_g}(G)=5.$   
As a consequence, ${\rm col_g}(\mathcal{P}_k) = 5$ for $k =7,8.$


\section{Essential Tools} 
For a graph $G,$ let $\Pi (G)$ be the set of linear orderings of $V(G).$ 
The digraph $G_L$ with respect to $L\in \Pi(G)$ is obtained from $G$ by 
orienting an edge $uv$ in $G$ with $u >_L v$ into an arc $(u,v).$ 

For a vertex $u,$ we denote the set of neighbors of $u$ in $G$ by $N_G(u),$ 
the set of out-neighbors of $u$ in $G_L$ by $N^+_{G_L}(u),$ 
and the set of in-neighbors of $u$ in $G_L$ by $N^-_{G_L}(u).$ 
Let $d_G(u):= |N_G(u)|, d^+_{G_L}(u):= |N^+_{G_L}(u)|,$ and 
$d^-_{G_L}(u):= |N^-_{G_L}(u)|.$ 

We define $ V^+_{G_L}(u) := \{v \in V(G_L): v <_L u\}$ 
and $ V^-_{G_L}(u) := \{v \in V(G_L): v >_L u\}.$  
Moreover,  $N^+_{G_L}[u] := N^+_{G_L}(u) \cup \{u\}, 
N^-_{G_L}[u] := N^-_{G_L}(u) \cup \{u\}, 
V^+_{G_L}[u] := V^+_{G_L}(u) \cup \{u\},$ and  
$V^-_{G_L}[u] := V^-_{G_L}(u) \cup \{u\}.$ 
The subscripts maybe omitted if $G$ or $G_L$ is clear from the context. 

In \cite{K00}, Kierstead introduces the activation strategy which gives the bound for $\cn(G)$ 
in terms of the following parameters. 
 
\begin{Def}\label{d1} 
For a given graph $G,$ we say that $M$ is a \textit{matching from} $A$ \textit{to} $B$ 
if $M$ covers all vertices in A and each edge in $M$ joins a vertex in $A$ and 
a vertex in $B - A.$ 

For a vertex $u$ in $G_L,$ we define $m(u, L,G)$ to be the size of a largest 
$Z$  such that $Z\subseteq N^{-}[u] $ with a partition $Z = X \cup Y$ and there exist matchings 
$M$ from $X \subseteq  N^{-}[u] $ to $V^+(u)$ 
and $N$ from $Y \subseteq  N^{-}(u) $ to $V^+[u].$  

Let 
$$ r(u,L,G) := d^+_{G_L}(u) + m(u,L,G),  $$
$$  r(L,G) := \max_{u \in V(G)} r(u,L,G), $$
$$  r(G) := \min_{L \in \Pi(G)} r(L,G).$$ 
\end{Def} 

\begin{Th} (Kierstead \cite{K00})
$\cn(G) \leq 1+r(G).$
\end{Th}

\section{Upper bounds for the game coloring number of planar graphs with girth at least 7} 
\begin{Th}\label{T2} 
If $G$ is a planar graph with girth at least $7,$ then $\cn(G) \leq 5.$ 
\end{Th} 
\begin{Pf}
Fix a planar graph embedding of $G.$  
It suffices to construct a linear ordering $L$ such that $r(L,(G)) \leq 4$ as follows. 
Initially, we have a set of chosen vertices $C := \emptyset$ 
and a set of unchosen vertices $U := V(G).$ 
At each stage of the construction, we choose a vertex $u \in U$, 
and replace $U$ by $U = \{u\}$ and $C$ by $C= C \cup \{u\}.$ 
Define a linear order $L$ by $u <_L v$ if we choose $v$ before $u.$     

It is well known that for a planar graph $G$ with girth at least $4,$ there exists a 
vertex $u$ with degree at most $3.$ 
If $C$ is empty, then choose a vertex of degree at most $3$ as $u.$ 
Suppose that $C$ is not empty. 
Let $H$ be the graph obtained from $G$ by \\ 
(i) deleting all edges between vertices in $C,$ 
(ii) deleting each vertex $x \in C$ such that $|N_G(x) \cap U| \leq 3,$  
(iii) if $x\in C$ and $|N_G(x) \cap U| =2,$ then we add an edge between 
two vertices in $N_G(x) \cap U,$   
(iv) if $x\in C$ and $|N_G(x) \cap U| =3,$ then we add two new edges between 
three vertices in $N_G(x) \cap U$ to form a path.  
Clearly, $H$ is a planar graph. 
Note that by $G$ has girth at least $4,$ a new edge in $H$ joins two vertices that are not adjacent in $G.$ 
Since the girth of $G$ is at least $7,$ it follows 
that the girth of $H$ is at least $4.$ 

Note that each $v \in C$ has $d_H(v) \geq 4.$ 
Thus there is a vertex $u \in U$ with $d_H(u) \leq 3.$ 
Choose such $u.$ 

Let   
$$S := \{y \in U: uy \in E(G) \},$$  
$$S' :=\{y \in U: uy \in E(H) - E(G) \},$$ 
$$A := \{x \in C: ux\in E(H) \textrm{ and } d_H(x) \geq 4\}.$$ 
Let $\sigma , \sigma',$ and $\alpha$ 
be the cardinalities of $S, S',$ and $A,$ respectively. 

Then 
$$ d(u) = \sigma + \sigma' + \alpha \leq 3.$$ 

Consider a set $Z \subseteq N^{-}[u] $ with $|Z| = m(u,L,G)$ 
such that there exists a partition $Z = X \cup Y$ and there exist matchings 
$M$ from $X \subseteq  N^{-}[u] $ to $V^+(u)$ 
and $N$ from $Y \subseteq  N^{-}(u) $ to $V^+[u].$  
Let $Z_1:= Z\cap A, Z_2 = Z \cap \{x \in C: d_H(x) \leq 3\},$ and $Z_3= Z \cap \{a\}.$ 
Then the sets $Z_i$ for $i=1,2,3$ are mutually disjoint 
but some sets maybe empty.   

Clearly, $d^+(u)=\sigma.$ 
We show that $m(u,L,G) \leq  \alpha + \sigma' +|Z_4|.$ 
From definitions, $|Z_1| \leq \alpha.$ 
It remains to show that $|Z_2| \leq \sigma'.$ 
Let $N' = N - \{uy \in N\}.$ 
Suppose $bx \in M$ and $by \in N'$ where $d_H(x) \leq 3$ 
and  $d_H(y) \leq 3.$  
Then $uxbyu$ is a cycle of length $4$ in $G,$ a contradiction. 
It follows that $b$ where $ub \in S'$ is the endpoint of at most one 
edge in $M\cup N'.$ 
Moreover, one endpoint of each edge in $M \cup N'$ is 
an endpoint of an edge in $S'.$ 
Thus $|Z_2| \leq \sigma'.$  
Now $r(u,L,G) = d^+(u) + m(u,L,G) \leq \sigma + \alpha + \sigma' +|Z_4|
\leq \sigma + \alpha + \sigma' +1.$      
Since  $\sigma+\sigma'+ \alpha \leq d_H(u)\leq 3,$   we have $r(u,L,G) \leq 4.$ 
Hence $\cn(G) \leq 5.$ 
\end{Pf}

\section{Lower bounds for the game coloring number of planar graphs with girth at most 8} 

\begin{Th}\label{T3} 
For $3 \leq k \leq 8,$ there is a planar graph $G$ with girth $k$ such that $\cn(G)= 5.$ 
As a consequence, $\cn(\mathcal{P}_k) \geq 5.$ 
\end{Th} 
\begin{Pf}
Let $H_1$ be a hexagon (a cycle of length 6). 
We construct $H_n$ ($n\geq 2$) from $H_{n-1}$ by adding 
$6(n-1)$ hexagon around $H_{n-1}.$ 

\begin{figure} [ht]
	\begin{center}
\begin{pspicture}(0,-1.87)(9.26,1.87)
\psline[linewidth=0.04cm](0.43,0.6)(1.21,0.6)
\psline[linewidth=0.04cm](0.43,0.6)(0.05,0.0)
\psline[linewidth=0.04cm](1.21,0.62)(1.61,0.02)
\psline[linewidth=0.04cm](0.05,0.04)(0.39,-0.58)
\psline[linewidth=0.04cm](1.59,0.02)(1.27,-0.58)
\psline[linewidth=0.04cm](0.39,-0.58)(1.29,-0.58)
\psdots[dotsize=0.14](0.45,0.58)
\psdots[dotsize=0.14](1.23,0.58)
\psdots[dotsize=0.14](0.07,0.0)
\psdots[dotsize=0.14](1.59,0.02)
\psdots[dotsize=0.14](0.39,-0.58)
\psdots[dotsize=0.14](1.29,-0.58)
\psline[linewidth=0.04cm](6.75,1.8)(7.53,1.8)
\psline[linewidth=0.04cm](6.75,1.8)(6.37,1.2)
\psline[linewidth=0.04cm](7.53,1.82)(7.93,1.22)
\psline[linewidth=0.04cm](6.37,1.24)(6.71,0.62)
\psline[linewidth=0.04cm](7.91,1.22)(7.59,0.62)
\psdots[dotsize=0.14](6.77,1.78)
\psdots[dotsize=0.14](7.55,1.78)
\psdots[dotsize=0.14](6.39,1.2)
\psdots[dotsize=0.14](7.91,1.22)
\psline[linewidth=0.04cm](6.75,0.62)(7.53,0.62)
\psline[linewidth=0.04cm](6.75,0.62)(6.37,0.02)
\psline[linewidth=0.04cm](7.53,0.64)(7.93,0.04)
\psline[linewidth=0.04cm](6.37,0.06)(6.71,-0.56)
\psline[linewidth=0.04cm](7.91,0.04)(7.59,-0.56)
\psdots[dotsize=0.14](6.77,0.6)
\psdots[dotsize=0.14](7.55,0.6)
\psdots[dotsize=0.14](6.39,0.02)
\psdots[dotsize=0.14](7.91,0.04)
\psline[linewidth=0.04cm](6.75,-0.58)(7.53,-0.58)
\psline[linewidth=0.04cm](6.75,-0.58)(6.37,-1.18)
\psline[linewidth=0.04cm](7.53,-0.56)(7.93,-1.16)
\psline[linewidth=0.04cm](6.37,-1.14)(6.71,-1.76)
\psline[linewidth=0.04cm](7.91,-1.16)(7.59,-1.76)
\psdots[dotsize=0.14](6.77,-0.6)
\psdots[dotsize=0.14](7.55,-0.6)
\psdots[dotsize=0.14](6.39,-1.18)
\psdots[dotsize=0.14](7.91,-1.16)
\psline[linewidth=0.04cm](6.71,-1.78)(7.61,-1.78)
\psdots[dotsize=0.14](6.71,-1.78)
\psdots[dotsize=0.14](7.61,-1.78)
\psline[linewidth=0.04cm](5.55,1.22)(6.33,1.22)
\psline[linewidth=0.04cm](5.55,1.22)(5.17,0.62)
\psline[linewidth=0.04cm](5.17,0.66)(5.51,0.04)
\psline[linewidth=0.04cm](5.51,0.04)(6.41,0.04)
\psdots[dotsize=0.14](5.57,1.2)
\psdots[dotsize=0.14](5.19,0.62)
\psdots[dotsize=0.14](5.51,0.04)
\psline[linewidth=0.04cm](5.51,0.02)(5.13,-0.58)
\psline[linewidth=0.04cm](5.13,-0.54)(5.47,-1.16)
\psline[linewidth=0.04cm](5.47,-1.16)(6.37,-1.16)
\psdots[dotsize=0.14](5.15,-0.58)
\psdots[dotsize=0.14](5.47,-1.16)
\psline[linewidth=0.04cm](7.97,1.24)(8.75,1.24)
\psline[linewidth=0.04cm](8.75,1.26)(9.15,0.66)
\psline[linewidth=0.04cm](9.13,0.66)(8.81,0.06)
\psline[linewidth=0.04cm](7.93,0.06)(8.83,0.06)
\psdots[dotsize=0.14](8.77,1.22)
\psdots[dotsize=0.14](9.13,0.66)
\psdots[dotsize=0.14](8.83,0.06)
\psline[linewidth=0.04cm](8.79,0.06)(9.19,-0.54)
\psline[linewidth=0.04cm](9.17,-0.54)(8.85,-1.14)
\psline[linewidth=0.04cm](7.97,-1.14)(8.87,-1.14)
\psdots[dotsize=0.14](9.17,-0.54)
\psdots[dotsize=0.14](8.87,-1.14)
\end{pspicture} 	
	\caption{$H_1$ and $H_2$}
	\end{center}
\end{figure}
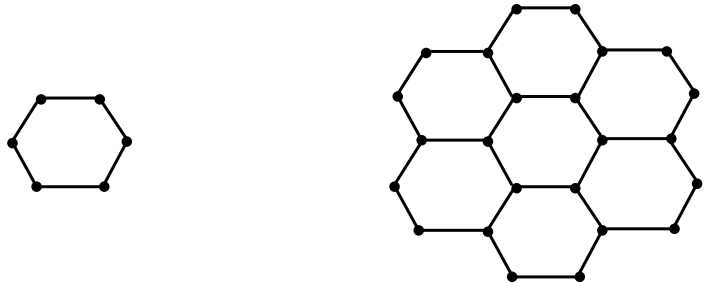

We construct $G_n$ from $H_n$ as follows. 
For each $v \in V(H_n),$ we add a vertex $a_v$ and an edge between $a_v$ and $v.$ 
Furthermore, we subdivide each horizontal edge of $H_n.$ 
The resulting graph $G_n$ is a planar graph with girth $8.$ 

\begin{figure} [ht]
	\begin{center}
\begin{pspicture}(0,-1.88)(4.18,1.88)
\psline[linewidth=0.04cm](1.67,1.81)(2.45,1.81)
\psline[linewidth=0.04cm](1.67,1.81)(1.29,1.21)
\psline[linewidth=0.04cm](2.45,1.83)(2.85,1.23)
\psline[linewidth=0.04cm](1.29,1.25)(1.63,0.63)
\psline[linewidth=0.04cm](2.83,1.23)(2.51,0.63)
\psdots[dotsize=0.14](1.69,1.79)
\psdots[dotsize=0.14](2.47,1.79)
\psdots[dotsize=0.14](1.31,1.21)
\psdots[dotsize=0.14](2.83,1.23)
\psline[linewidth=0.04cm](1.67,0.63)(2.45,0.63)
\psline[linewidth=0.04cm](1.67,0.63)(1.29,0.03)
\psline[linewidth=0.04cm](2.45,0.65)(2.85,0.05)
\psline[linewidth=0.04cm](1.29,0.07)(1.63,-0.55)
\psline[linewidth=0.04cm](2.83,0.05)(2.51,-0.55)
\psdots[dotsize=0.14](1.69,0.61)
\psdots[dotsize=0.14](2.47,0.61)
\psdots[dotsize=0.14](1.31,0.03)
\psdots[dotsize=0.14](2.83,0.05)
\psline[linewidth=0.04cm](1.67,-0.57)(2.45,-0.57)
\psline[linewidth=0.04cm](1.67,-0.57)(1.29,-1.17)
\psline[linewidth=0.04cm](2.45,-0.55)(2.85,-1.15)
\psline[linewidth=0.04cm](1.29,-1.13)(1.63,-1.75)
\psline[linewidth=0.04cm](2.83,-1.15)(2.51,-1.75)
\psdots[dotsize=0.14](1.69,-0.59)
\psdots[dotsize=0.14](2.47,-0.59)
\psdots[dotsize=0.14](1.31,-1.17)
\psdots[dotsize=0.14](2.83,-1.15)
\psline[linewidth=0.04cm](1.63,-1.77)(2.53,-1.77)
\psdots[dotsize=0.14](1.63,-1.77)
\psdots[dotsize=0.14](2.53,-1.77)
\psline[linewidth=0.04cm](0.47,1.23)(1.25,1.23)
\psline[linewidth=0.04cm](0.47,1.23)(0.09,0.63)
\psline[linewidth=0.04cm](0.09,0.67)(0.43,0.05)
\psline[linewidth=0.04cm](0.43,0.05)(1.33,0.05)
\psdots[dotsize=0.14](0.49,1.21)
\psdots[dotsize=0.14](0.11,0.63)
\psdots[dotsize=0.14](0.43,0.05)
\psline[linewidth=0.04cm](0.43,0.03)(0.05,-0.57)
\psline[linewidth=0.04cm](0.05,-0.53)(0.39,-1.15)
\psline[linewidth=0.04cm](0.39,-1.15)(1.29,-1.15)
\psdots[dotsize=0.14](0.07,-0.57)
\psdots[dotsize=0.14](0.39,-1.15)
\psline[linewidth=0.04cm](2.89,1.25)(3.67,1.25)
\psline[linewidth=0.04cm](3.67,1.27)(4.07,0.67)
\psline[linewidth=0.04cm](4.05,0.67)(3.73,0.07)
\psline[linewidth=0.04cm](2.85,0.07)(3.75,0.07)
\psdots[dotsize=0.14](3.69,1.23)
\psdots[dotsize=0.14](4.05,0.67)
\psdots[dotsize=0.14](3.75,0.07)
\psline[linewidth=0.04cm](3.71,0.07)(4.11,-0.53)
\psline[linewidth=0.04cm](4.09,-0.53)(3.77,-1.13)
\psline[linewidth=0.04cm](2.89,-1.13)(3.79,-1.13)
\psdots[dotsize=0.14](4.09,-0.53)
\psdots[dotsize=0.14](3.79,-1.13)
\psdots[dotsize=0.14](2.09,1.79)
\psdots[dotsize=0.14](2.07,0.63)
\psdots[dotsize=0.14](2.07,-0.59)
\psdots[dotsize=0.14](2.05,-1.79)
\psdots[dotsize=0.14](0.87,1.23)
\psdots[dotsize=0.14](3.27,1.25)
\psdots[dotsize=0.14](3.27,0.07)
\psdots[dotsize=0.14](0.89,0.05)
\psdots[dotsize=0.14](0.83,-1.15)
\psdots[dotsize=0.14](3.29,-1.15)
\end{pspicture} 
	\caption{$G_2$ (we omit the vertices of $A$)}
	\end{center}
\end{figure}
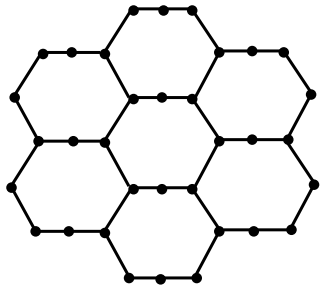

Let \\
$$A := \{a_v \in V(H_n): v \in V(H_n)\},$$ 
$$B := V(G_n) - (V(H_n) \cup A),$$  
$$V_I := \{v \in V(H_n): v \textrm{ is not incident to the unbounded face of } H_n\},$$
$$V_O := \{v \in V(H_n): v \textrm{ is incident to the unbounded face of } H_n\}.$$

Then we have \\
$$|B| = 3n^2-n,$$  
$$|V_I| = 6n^2-12 n+ 6,$$ 
$$|V_O| = 12 n- 6.$$ 

First, Bob marks all the vertices of $B$ and $V_O.$  
Note that $|B|+|V_O| +1 < |V_I|$ for $n \geq 9,$ 
It follows that when all the vertices of $B$ and $V_O$ are marked,  
there are at least two vertices of $V_I$ are unmarked. 
Consider Bob's turn immediately after all the vertices of $B \cup V_O$ are marked. 
Suppose $v_1,v_2,\ldots, v_k$ are all unmarked vertices in $V_I.$  
Bob can force all $a_{v_1},a_{v_2},\ldots,a_{v_k}$ to be marked before 
all of $v_1,v_2,\ldots, v_k$ are marked. 
Thus the last unmarked vertex, say $v$, in $V_I$ satisfies $b(v)+1 = 5.$ 
Therefore $\cn(G_n) \geq 5.$ 
Since $\Delta(G_n)=4,$ we have $\cn(G_n) = 5.$ 

For $k \leq 8,$ one can see that  $\cn(G_n \cup C_k)=5.$ 
It follows that $\cn(\mathcal{P}_k) \geq 5.$  
\end{Pf}

\begin{cor}\label{cor1} 
$\cn(\mathcal{G}_k) = 5$ for $k= 7,8.$  
\end{cor} 
\begin{Pf}
From Theorems \ref{T2}, \ref{T3}, and the fact that 
$\cn(\mathcal{P}_8) \leq 5$ \cite{WZ11},  Corollary \ref{cor1} follows.  
\end{Pf}

\section*{Acknowledgments} 
This work draws substantially on ideas from the paper by Kierstead \cite{K00} 
and the paper by Sekiguchi \cite{S14}.   
Without those ideas, this work would not be possible. 
We would like to thank them for their papers.  

The second author was supported by the Commission 
on Higher Education and the Thailand Research Fund 
under grant RSA5780014.


\begin{thebibliography}{99}

\bibitem{B91}
H.L. Bodlaender,  
On the complexity of some colouring games, 
\emph{J. Found. Comput. Sci.} 2(1991), 133--147.

\bibitem{BIKS09} 
O.V. Borodin, A.O. Ivanova, A.V. Kostochka, N.N. Sheikh, 
Decomposing of quadrangle-free planar graphs, 
\emph{Discuss. Math. Graph Theory} 29(2009), 87--99.

\bibitem{FKKT93} 
U. Faigle, U. Kern, H.A. Kierstead, W.T. Trotter, 
The game chromatic number of some classes of graphs graphs, 
\emph{Ars Combin.} 35(1993) 143--150 

\bibitem{GZ99} 
D. Guan, X. Zhu, 
The game chromatic number of outerplanar graphs, 
\emph{J. Graph Theory} 30(1999), 67--70.

\bibitem{HEtal02} 
W. He, X. Hou, K. Lih, J. Shao, W. Wang, X. Zhu, 
Edge-partitions of planar graphs and their game coloring numbers, 
\emph{J. Graph Theory} 41(2002), 307--317.

\bibitem{K06}
D. Kleitman, 
Partitioning the edges of girth 6 planar graph into those of a forest and those of a set of disjoint paths and cycles, manuscript, 2006.


\bibitem{K00}
H. A. Kierstead, 
A simple competitive graph colouring algorithm, 
\emph{J. Combin. Theory Ser. B} 78(2000), 57--68. 

\bibitem{KY05}
H. A. Kierstead, D. Yang
Very asymmetric marking games, 
\emph{Order} 22(2005), 93--107. 

\bibitem{S14}
Y. Sekiguchi, 
The game coloring number of planar graphs with given girth, 
\emph{Discrete Math.} 330(2014),  11--16.

\bibitem{WZ11}
Y. Wang, Q. Zhang,  
Decomposing a planar graph with girth at least 8 into a forest and a matching, 
\emph{Discrete Math.} 311(2011),  844--849. 


\bibitem{WZ08}
J. Wu, X. Zhu, 
Lower bounds for the game colouring number of partial $k$-trees and planar graphs, 
\emph{Discrete Math.} 308(2008),  2637--2642. 

\bibitem{Z99}
X. Zhu, 
The game colouring number of planar graphs, 
\emph{J. Combin. Theory Ser. B} 75(1999), 245--258.

\bibitem{Z00}
X. Zhu, 
Game colouring number of pseudo partial $k$-trees and planar graphs, 
\emph{Discrete Math.} 215(2000),  245--262. 

\bibitem{Z08}
X. Zhu, 
Refined activation strategy for the marking game, 
\emph{J. Combin. Theory Ser. B} 98(2008), 1--18.


\end{thebibliography}
\end{document}